\documentclass{article}
\usepackage{amssymb}
\usepackage{amsfonts}
\usepackage{amsmath}

\newtheorem{theorem}{Theorem}

\newtheorem{corollary}[theorem]{Corollary}

\newtheorem{definition}[theorem]{Definition}
\newtheorem{example}[theorem]{Example}

\newtheorem{lemma}[theorem]{Lemma}
\newtheorem{notation}[theorem]{Notation}

\newtheorem{proposition}[theorem]{Proposition}
\newtheorem{remark}[theorem]{Remark}

\newenvironment{proof}[1][Proof]{\noindent\textbf{#1.} }{\ \rule{0.5em}{0.5em}}
\input{tcilatex}

\begin{document}

\title{Euler Estimates for Rough Differential Equations}
\author{Peter Friz \and Nicolas Victoir}
\maketitle

\begin{abstract}
We consider controlled differential equations and give new estimates for
higher order Euler schemes. Our proofs are inspired by recent work of A. M.
Davie who considers first and second order schemes. In order to implement
the general case we make systematic use of geodesic approximations in the
free nilpotent group.

As application, we can control moments of solutions to rough path
differential equations (RDEs) driven by random rough paths with sufficient
integrability and have a criteria for $L^{q}$-convergence in the Universal
Limit Theorem. We also obtain Azencott type estimates and asymptotic
expansions for random RDE solution. \ When specialized to RDEs driven by
Enhanced Brownian motion, we (mildly) improve classic estimates for
diffusions in the small time limit.
\end{abstract}

\section*{Introduction}

\bigskip We consider controlled differential equations of the form%
\begin{equation*}
dy=V_{1}\left( y\right) dx^{1}+...+V_{d}\left( y\right) dx^{d}
\end{equation*}%
on the time interval $\left[ 0,1\right] $. When the $\mathbb{R}^{d}$-valued
driving signal $x$ and the vector fields are Lipschitz continuous then there
exists a unique solution for every starting point $y_{0}$. Building on
recent work by A. M. Davie \cite{Da}, we control the H\"{o}lder norm of $y$
in terms of suited H\"{o}lder norms of $x$ and its iterated integrals and
also obtain higher order Euler estimates. To be able to go beyond the first
and second order case discussed by Davie we use ideas from sub-Riemannian
geometry. En passant, we construct \textit{geodesic approximations} which
respect the geometry of the vector fields. Finally, by passing to the limit,
all estimates extend to solutions of rough path differential equations
(RDEs) in the sense of T. Lyons.

As application, we show that random\ RDE solutions driven by a sufficiently
integrable geometric rough path are in $L^{q}\left( \Omega \right) $ for all 
$q<\infty $. The examples we have in mind are RDEs\ driven by \textit{%
Enhanced Brownian motion} (in which case we are effectively dealing with a
Stratonovich stochastic differential equation) as well as certain \textit{%
Enhanced Gaussian -} and \textit{Enhanced Markov processes.}

When specialized to RDEs driven by Enhanced Brownian motion, the above
mentioned Euler estimates are closely related to classic estimates by
Azencott. In fact, at least in absence of a drift vector field, our main
Euler estimate sharpens a key result in \cite{Az} when applied to \textit{%
diffusions en temps petit,} which is precisely the result which is quoted
and used in Ben Arous' and, later, Castell's work \cite{BA, Ca}.

As the reader may suspect, the robustness of the rough path approach allows
to obtain Azencott type estimates for \textit{arbitrary} random\ RDE
solutions driven by a sufficiently integrable geometric rough path. The case
of RDEs driven by EBM, as discussed in the last paragraph, is just a special
case of this general result.

\begin{notation}
The dimensions of $\mathbb{R}^{d},\mathbb{R}^{e}$ are fixed and will not
appear explicitly when we write out the dependence of constants. In general,
constants that appear in lemmas, propositions, theorems etc have an index
that matches the number of the statement. In the proofs we indicate changing
constants by a running upper index.
\end{notation}

\section{Preliminaries: Euler scheme of order $N$}

\bigskip Define $T^{\left( N\right) }=\oplus _{k=0}^{N}\left( \mathbb{R}%
^{d}\right) ^{\otimes k}$, by convention $\left( \mathbb{R}^{d}\right)
^{\otimes 0}\equiv \mathbb{R}$. Let $x$ be an $\mathbb{R}^{d}$-valued
Lipschitz path and define the $k^{th}$ iterated integrals of the path
segment $x|_{\left[ s,t\right] }$ as%
\begin{equation*}
\mathbf{g}^{k,i_{1},\cdots
,i_{k}}:=\int_{s}^{t}\int_{s}^{u_{k}}...%
\int_{s}^{u_{2}}dx_{u_{1}}^{i_{1}}...dx_{u_{k}}^{ik}.
\end{equation*}

and so that $\mathbf{g}^{k}=\left( \mathbf{g}^{k,i_{1},\cdots ,i_{k}}\right)
_{i_{1},\cdots ,i_{k}\in \left\{ 1,...,d\right\} }\in $ $\left( \mathbb{R}%
^{d}\right) ^{\otimes k}$. For later convenience set $\mathbf{g}^{0}=1\in
\left( \mathbb{R}^{d}\right) ^{\otimes 0}\equiv \mathbb{R}$. We then define
the (step-$N$) signature of the path segment $x|_{\left[ s,t\right] }$ as 
\begin{equation*}
\mathbf{x}_{s,t}\equiv S_{N}\left( x\right) _{s,t}\equiv 1+\sum_{k=1}^{N}%
\mathbf{g}^{k}\in T^{\left( N\right) }\left( \mathbb{R}^{d}\right) .
\end{equation*}

\bigskip

We say that a vector field is in $\mathrm{Lip}^{\gamma }\left( \mathbb{R}%
^{e}\right) $ if it has $\lfloor \gamma \rfloor $ bounded derivatives and
the $\lfloor \gamma \rfloor ^{th}$-derivative is $\left\{ \gamma \right\} $-H%
\"{o}lder continuous.

\begin{definition}
Given vector fields $V_{1},...,V_{d}\in \mathrm{Lip}^{1}\left( \mathbb{R}%
^{e}\right) $ (= bounded \& Lipschitz continuous vector fields) and a $%
\mathbb{R}^{d}$-valued Lipschitz path $x$ on $\left[ 0,1\right] ,$ we let $%
y=\pi \left( 0,y_{0};x\right) $ denote the unique solution to the (control)\
ODE%
\begin{equation*}
dy_{t}=\sum_{i=1}^{d}V_{i}\left( y_{t}\right) dx_{t}^{i}\equiv V\left(
y_{t}\right) dx_{t},\,\,\,t\in \left[ 0,1\right] .
\end{equation*}%
started at $y_{0}.$
\end{definition}

The following lemma is left as a simple exercise.

\begin{lemma}
\label{boundY}Assume that $\left( V_{i}\right) _{1\leq i\leq d}\in \mathrm{%
Lip}^{1}\left( \mathbb{R}^{e}\right) .$ Let $x$ be an $\mathbb{R}^{d}$%
-valued Lipschitz path on $\left[ 0,1\right] $ and let $y_{t}=\pi \left(
0,y_{0};x\right) _{t}$. Then, for all $0\leq s\leq t\leq 1,$%
\begin{equation*}
\left\vert y_{s,t}\right\vert \leq C_{\ref{boundY}}\int_{s}^{t}\left\vert
dx_{r}\right\vert
\end{equation*}%
where $C_{\ref{boundY}}$ depends on (the Lipschitz norm of) the vector
fields $V_{1},...V_{d}$.
\end{lemma}

Let us now define the Euler approximation of order $N$ to a control ODE of
the above type. To this end, let $H$ denote the identity function on $%
\mathbb{R}^{e}$ and recall the identification of vector fields with first
order differential operators.

\begin{definition}
Given $\left( V_{i}\right) _{1\leq i\leq d}$ $\in \mathrm{Lip}^{N}\left( 
\mathbb{R}^{e}\right) ,$ $\mathbf{g}\in T^{\left( N\right) }\left( \mathbb{R}%
^{d}\right) $ and $y\in \mathbb{R}^{e}$ we call%
\begin{equation*}
I^{y,N,\mathbf{g}}:=I\left[ y,N,\mathbf{g}\right] :=\sum_{k=1}^{N}\sum 
_{\substack{ i_{1},...,i_{k}  \\ \in \left\{ 1,...,d\right\} }}%
V_{i_{1}}\cdots V_{i_{k}}H\left( y\right) \mathbf{g}^{k,i_{1},\cdots ,i_{k}},
\end{equation*}%
the (increment of) the step-$N$ Euler scheme.
\end{definition}

This definition is explained by

\begin{lemma}
\label{Euler}Assume that $\left( V_{i}\right) _{1\leq i\leq d}\in \mathrm{Lip%
}^{N}\left( \mathbb{R}^{e}\right) .$ Let $x$ be an $\mathbb{R}^{d}$-valued
Lipschitz path on $\left[ 0,1\right] $ and let $y_{t}=\pi \left(
0,y_{0};x\right) _{t}$. Then, for all $0\leq s\leq t\leq 1,$%
\begin{eqnarray*}
&&y_{s,t}-I^{y_{s},N,S_{N}\left( x\right) _{s,t}} \\
&=&\sum_{\substack{ i_{1},...,i_{N}  \\ \in \left\{ 1,...,d\right\} }}%
\int_{s<r_{1}<...<r_{N}<t}\left[ V_{i_{1}}\cdots V_{i_{N}}H\left(
y_{r_{1}}\right) -V_{i_{1}}\cdots V_{i_{N}}H\left( y_{s}\right) \right]
dx_{r_{1}}^{i_{1}}\cdots dx_{r_{N}}^{i_{N}}
\end{eqnarray*}%
and there exists a constant $C_{\ref{Euler}}$ depending on $N$ and $%
V_{1},...,V_{d}$ such that%
\begin{equation*}
\left\vert y_{s,t}-I^{y_{s},N,\mathbf{x}_{s,t}}\right\vert \leq C_{\ref%
{Euler}}\left( \int_{s}^{t}\left\vert dx_{r}\right\vert \right) ^{N+1}.
\end{equation*}
\end{lemma}

\begin{proof}
Let $f$ be smooth and note that $V_{i}\in \mathrm{Lip}^{N}$ implies $%
V_{i_{1}}\cdots V_{i_{k}}f$ is $C^{1}$ for $1\leq k\leq N$. Iterated use of
the fundamental theorem of calculus gives 
\begin{eqnarray*}
f\left( y_{t}\right) &=&f\left( y_{s}\right) +\sum_{k=1}^{N-1}\sum 
_{\substack{ i_{1},...,i_{k}  \\ \in \left\{ 1,...,d\right\} }}%
\int_{s<r_{1}<...<r_{k}<t}V_{i_{1}}\cdots V_{i_{k}}f\left( y_{s}\right)
dx_{r_{1}}^{i_{1}}\cdots dx_{r_{k}}^{i_{k}} \\
&&+\sum_{\substack{ i_{1},...,i_{k}  \\ \in \left\{ 1,...,d\right\} }}%
\int_{s<r_{1}<...<r_{N}<t}V_{i_{1}}\cdots V_{i_{N}}f\left( y_{r_{1}}\right)
dx_{r_{1}}^{i_{1}}\cdots dx_{r_{N}}^{i_{N}}.
\end{eqnarray*}%
This first part is then proved by specializing to $f=H$. For the second
statement, lemma \ref{boundY} gives%
\begin{equation*}
\left\vert y_{s,t}\right\vert \leq C_{\ref{Euler}}^{1}\int_{s}^{t}\left\vert
dx_{r}\right\vert .
\end{equation*}

$\mathrm{Lip}^{N}$-regularity of the vector fields implies that $%
V_{i_{1}}..V_{i_{N}}H\left( \cdot \right) $ is Lipschitz and hence, for $%
r\in \left[ s,t\right] ,$ 
\begin{equation*}
\left\vert V_{i_{1}}..V_{i_{N}}H\left( y_{r}\right)
-V_{i_{1}}..V_{i_{N}}H\left( y_{s}\right) \right\vert \leq C_{\ref{Euler}%
}^{2}\int_{s}^{t}\left\vert dx_{r}\right\vert .
\end{equation*}%
This leads to%
\begin{equation*}
\left\vert \int_{s<r_{1}<...<r_{N}<t}\left[ V_{i_{1}}..V_{i_{N}}H\left(
y_{r_{1}}\right) -V_{i_{1}}..V_{i_{N}}H\left( y_{s}\right) \right]
dx_{r_{1}}^{i_{1}}\cdots dx_{r_{N}}^{i_{N}}\right\vert \leq C_{\ref{Euler}%
}^{3}\left( \int_{s}^{t}\left\vert dx_{r}\right\vert \right) ^{N+1}
\end{equation*}%
and summing over the indices finishes the proof.
\end{proof}

\section{Preliminaries II:\ Algebra of Iterated Integrals}

The set $T_{1}^{N}\left( \mathbb{R}^{d}\right) \equiv \left\{ \mathbf{g}\in
T^{N}\left( \mathbb{R}^{d}\right) :\mathbf{g}^{0}=1\right\} $ is a group
under \textit{truncated tensor multiplication:} if $\mathbf{g}=1+\mathbf{g}%
^{1}+...+\mathbf{g}^{N}\equiv 1+$ $\mathbf{\tilde{g}}$ and similar for $%
\mathbf{h}$ then for $k=0,...,N$%
\begin{equation*}
\left( \mathbf{g}\otimes \mathbf{h}\right) ^{k}=\sum_{i=0}^{k}\mathbf{g}%
^{i}\otimes \mathbf{h}^{k-i}.
\end{equation*}%
The neutral element is $\mathbf{e}=1=1+0+...+0$ and the inverse is given by
the usual power series calculus%
\begin{equation*}
\left( 1+\mathbf{\tilde{g}}\right) ^{-1}=1-\mathbf{\tilde{g}+\tilde{g}}%
^{\otimes 2}-...
\end{equation*}%
For every $\lambda \in \mathbb{R}$, the \textit{dilatation} map $\delta
_{\lambda }$ is defined componenwise by $\mathbf{g}^{k}$ $\mapsto \lambda
^{k}\mathbf{g}^{k}$, \thinspace \thinspace $k=0,...,N$. 
\begin{equation*}
\delta _{\lambda }:\left( \mathbf{g}^{k}\right) \mapsto \left( \lambda ^{k}%
\mathbf{g}^{k}\right) ,\text{ \ \ \ }\lambda \in \mathbb{R}.
\end{equation*}%
Obviously, $T_{1}^{N}\left( \mathbb{R}^{d}\right) $ is a Lie group. Its Lie
algebra can be identified with%
\begin{equation*}
T_{0}^{N}\left( \mathbb{R}^{d}\right) \equiv \left\{ \mathbf{\tilde{g}}\in
T^{N}\left( \mathbb{R}^{d}\right) :\mathbf{\tilde{g}}^{0}=0\right\} \text{, }%
\left[ \mathbf{\tilde{g},\tilde{h}}\right] =\mathbf{\tilde{g}}\otimes 
\mathbf{\tilde{h}}-\mathbf{\tilde{h}}\otimes \mathbf{\tilde{g}}
\end{equation*}%
and the exponential map with $\exp :T_{0}^{N}\left( \mathbb{R}^{d}\right)
\rightarrow T_{1}^{N}\left( \mathbb{R}^{d}\right) $, $\mathbf{\tilde{g}%
\mapsto }1+\mathbf{\tilde{g}+}\frac{1}{2!}\mathbf{\tilde{g}}^{\otimes 2}+...$

We recall some well-known facts. See \cite{Ly, FV04, Mo, Bd} for further
references.

\begin{proposition}[Chen, \protect\cite{Ly}]
\label{productOfSigs}Let $x:\left[ 0,1\right] \rightarrow $ $\mathbb{R}^{d}$
be Lipschitz continuous with (step-$N$) signatures $\mathbf{x}%
_{s,t}=S_{N}\left( x\right) _{s,t}$. Then%
\begin{equation}
S_{N}\left( x\right) _{s,t}\otimes S_{N}\left( x\right) _{t,u}=S_{N}\left(
x\right) _{s,u}.  \label{productOfSigsFormula}
\end{equation}
\end{proposition}

We define $G^{N}\left( \mathbb{R}^{d}\right) \equiv \exp \left( L^{N}\left( 
\mathbb{R}^{d}\right) \right) $ where%
\begin{equation*}
L=L^{N}\left( \mathbb{R}^{d}\right) =\mathbb{R}^{d}\oplus \left[ \mathbb{R}%
^{d},\mathbb{R}^{d}\right] \oplus \left[ \mathbb{R}^{d},\left[ \mathbb{R}%
^{d},\mathbb{R}^{d}\right] \right] \oplus ...\subset T_{0}^{N}\left( \mathbb{%
R}^{d}\right) ,
\end{equation*}

$G^{N}\left( \mathbb{R}^{d}\right) $ is a Lie subgroup of $T_{1}^{N}\left( 
\mathbb{R}^{d}\right) $ with respect to $\otimes $-multiplication and known
as \textit{step-}$N$\textit{\ nilpotent free group over }$\mathbb{R}^{d}$.

\begin{theorem}[Chow, \protect\cite{Mo}]
For every $\mathbf{g}\in G^{N}\left( \mathbb{R}^{d}\right) $ there exists an 
$\mathbb{R}^{d}$-valued Lipschitz path $x$ such that $S_{N}\left( x\right)
_{0,1}=\mathbf{g}$. $\ $More precisely, $G$ is the group generated by $%
\left\{ \exp \left( v\right) :v\in \mathbb{R}^{d}\right\} $ so that every $%
\mathbf{g}\in G$ is the signature of a (finite number of) concatenation of
straight path segments.
\end{theorem}

\begin{theorem}[Geodesic Existence, \protect\cite{Mo}]
For every $\mathbf{g}\in G^{N}\left( \mathbb{R}^{d}\right) $,%
\begin{equation*}
\left\Vert \mathbf{g}\right\Vert :=\inf \left\{ \int_{0}^{1}\left\vert \dot{%
\gamma}_{t}\right\vert dt:\gamma :\left[ 0,1\right] \rightarrow \mathbb{R}%
^{d}\text{ Lipschitz continuous, }\gamma \left( 0\right) =0,\text{ }%
S_{N}\left( \gamma \right) _{0,1}=\mathbf{g}\text{ }\right\}
\end{equation*}%
is finite and achieved at some minimizing Lipschitz continuous path $\gamma
^{\ast }$, i.e.%
\begin{equation*}
\left\Vert \mathbf{g}\right\Vert =\int_{0}^{1}\left\vert \dot{\gamma}%
_{t}^{\ast }\right\vert dt\text{ and }S_{N}\left( \gamma ^{\ast }\right)
_{0,1}=\mathbf{g}.
\end{equation*}%
Moreover, by simple reparametrization, we can state that for every\thinspace
\thinspace $s,t\in \mathbb{R}$ with $s<t$ \ then exists a Lipschitz path $%
x^{s,t}:\left[ s,t\right] \rightarrow \mathbb{R}^{d}$ with signature $%
\mathbf{g}$ and length $\left\Vert \mathbf{g}\right\Vert $:%
\begin{equation*}
S_{N}\left( x^{s,t}\right) _{s,t}=\mathbf{g}\ \text{\ \ and \ }%
\int_{s}^{t}\left\vert dx^{s,t}\right\vert =\left\Vert \mathbf{g}\right\Vert
.
\end{equation*}
\end{theorem}

\begin{remark}
$G^{N}\left( \mathbb{R}^{d}\right) $ can be given a subriemannian structure
so that the path $t\in \lbrack 0,1]\mapsto $ $S_{N}\left( \gamma ^{\ast
}\right) _{0,t}$ is a subriemannian geodesic which connects the unit $e$
with $\mathbf{g}\in G^{N}\left( \mathbb{R}^{d}\right) $, see \cite{Bd, Mo}.
Thus, strictly speaking, $\gamma ^{\ast }$ is not a geodesics but the
projection of a geodesic.
\end{remark}

The geodesic existence theorem has useful consequences. If $\mathbf{g},%
\mathbf{h}\in G^{N}\left( \mathbb{R}^{d}\right) ~$then \newline
(i) $\left\Vert \mathbf{g}\right\Vert =0$ iff $\mathbf{g}=e$, (ii) symmetry: 
$\left\Vert \mathbf{g}\right\Vert =\left\Vert \mathbf{g}^{-1}\right\Vert $,
(iii) sub-additivity $\left\Vert \mathbf{g}\otimes \mathbf{h}\right\Vert
\leq \left\Vert \mathbf{g}\right\Vert +\left\Vert \mathbf{h}\right\Vert $
and (iv) homogenity $\left\Vert \delta _{\lambda }\mathbf{g}\right\Vert
=\left\vert \lambda \right\vert \,\left\Vert \mathbf{g}\right\Vert $ for all 
$\lambda \in \mathbb{R}$, hold true. In particular, $d\left( \mathbf{g},%
\mathbf{h}\right) :=\left\Vert \mathbf{g}^{-1}\otimes \mathbf{h}\right\Vert $
defines a left-invariant metric on $G^{N}\left( \mathbb{R}^{d}\right) $, the 
\textit{Carnot-Caratheodory metric}.

\begin{theorem}[\protect\cite{Mo}]
(a) The topology induced by Carnot-Caratheodory metric coincides with the
manifold topology of $G^{N}\left( \mathbb{R}^{d}\right) $ and the trace
topology as as subset of $T_{1}^{N}\left( \mathbb{R}^{d}\right) $.\newline
(b) The map $\mathbf{g}\mapsto \left\vert \left\vert \mathbf{g}\right\vert
\right\vert $ is continuous in this topology.\newline
(c) The space $G^{N}\left( \mathbb{R}^{d}\right) $ with metric $d$ is Polish.
\end{theorem}

\begin{proposition}[\protect\cite{Go}]
Let $\left\vert \left\vert \left\vert .\right\vert \right\vert \right\vert
_{i}$ $(i=1,2)$ be continuous homogenous norms on $G^{N}\left( \mathbb{R}%
^{d}\right) $, that is, norms that satisfies properties $\left( i\right) $
and $\left( iv\right) $ and such that $\mathbf{g}\mapsto \left\vert
\left\vert \left\vert \mathbf{g}\right\vert \right\vert \right\vert _{i}$ is
continuous w.r.t. $\tau $. Then there exists a constant $c\in \lbrack
1,\infty )$ such that $\left\vert \left\vert \left\vert .\right\vert
\right\vert \right\vert _{1}\sim \left\vert \left\vert \left\vert
.\right\vert \right\vert \right\vert _{2}$ by which we mean%
\begin{equation*}
\frac{1}{c}\left\vert \left\vert \left\vert .\right\vert \right\vert
\right\vert _{2}\leq \left\vert \left\vert \left\vert .\right\vert
\right\vert \right\vert _{1}\leq c\left\vert \left\vert \left\vert
.\right\vert \right\vert \right\vert _{2}.
\end{equation*}
\end{proposition}

For instance,%
\begin{equation*}
\left\vert \left\vert \left\vert \mathbf{g}\right\vert \right\vert
\right\vert \equiv \max_{k=1,...,N}\left\vert \mathbf{g}^{k}\right\vert
^{1/k}.
\end{equation*}
provides a useful example of a continuous homogenous norm on $G^{N}\left( 
\mathbb{R}^{d}\right) $ other than $\left\Vert \cdot \right\Vert $.

\section{Preliminaries III:\ Geometric (H\"{o}lder) Rough Paths}

Here, and in the remainder of this paper, we work exclusively with H\"{o}%
lder modulus\bigskip 
\begin{equation*}
\omega \left( s,t\right) \equiv t-s.
\end{equation*}
Let $p\in \lbrack 1,\infty )$. A path $\mathbf{x}$\textbf{\ }from $\left[ 0,1%
\right] $ to $G^{N}\left( \mathbb{R}^{d}\right) $ is \thinspace $1/p$-H\"{o}%
lder continuous if for all $s,t\in \left[ 0,1\right] $%
\begin{equation*}
\left\Vert \mathbf{x}_{s,t}\right\Vert \leq C\omega \left( s,t\right) ^{1/p}
\end{equation*}%
for some constant $C.$ This class is denoted by $C^{1/p\text{-H\"{o}lder}%
}\left( \left[ 0,1\right] ,G^{N}\left( \mathbb{R}^{d}\right) \right) $. We
can restrict attention to paths with pinned starting point. The (homogenous) 
$1/p$-H\"{o}lder "norm" (there is no linear space here) on $C^{1/p\text{-H%
\"{o}lder}}\left( \left[ 0,1\right] ,G^{N}\left( \mathbb{R}^{d}\right)
\right) $\ is defined by%
\begin{equation*}
\left\Vert \mathbf{x}\right\Vert _{1/p\text{-H\"{o}lder;}\left[ 0,1\right]
}=\left\Vert \mathbf{x}\right\Vert _{1/p\text{-H\"{o}lder}}=\sup_{0\leq
s<t\leq 1}\frac{\left\Vert \mathbf{x}_{s,t}\right\Vert }{\omega \left(
s,t\right) ^{1/p}}
\end{equation*}%
and there is a $1/p$-H\"{o}lder metric based on the CC-metric,%
\begin{equation*}
d_{1/p\text{-H\"{o}lder;}\left[ 0,1\right] }\left( \mathbf{x,\tilde{x}}%
\right) =d_{1/p\text{-H\"{o}lder}}\left( \mathbf{x,\tilde{x}}\right)
=\sup_{0\leq s<t\leq 1}\frac{d\left( \mathbf{x}_{s,t},\mathbf{\tilde{x}}%
_{s,t}\right) }{\omega \left( s,t\right) ^{1/p}}.
\end{equation*}%
We also set%
\begin{equation*}
d_{\infty ;\left[ 0,1\right] }\left( \mathbf{x,\tilde{x}}\right) =d_{\infty
}\left( \mathbf{x,\tilde{x}}\right) =\sup_{0\leq s<t\leq 1}d\left( \mathbf{x}%
_{s,t},\mathbf{\tilde{x}}_{s,t}\right) .
\end{equation*}

\begin{theorem}[\protect\cite{FV04}]
\label{ThmOnHoelderRoughPathsProperties}(i) $C^{1/p\text{-H\"{o}lder}}\left( %
\left[ 0,1\right] ,G^{N}\left( \mathbb{R}^{d}\right) \right) $ is a complete
metric space under the metric $d_{1/p\text{-H\"{o}lder}}$.\newline
(ii) Every $1/p$-H\"{o}lder continuous path $\mathbf{x}\in C^{1/p\text{-H%
\"{o}lder}}\left( \left[ 0,1\right] ,G^{N}\left( \mathbb{R}^{d}\right)
\right) $ can be approximated by Lipschitz paths $x_{n}:\left[ 0,1\right]
\rightarrow \mathbb{R}^{d}$ in the sense that%
\begin{equation*}
S_{N}\left( x_{n}\right) \rightarrow \mathbf{x}\text{ \ uniformly on }\left[
0,1\right]
\end{equation*}%
and $\sup_{n}\left\Vert S_{N}\left( x_{n}\right) \right\Vert _{1/p\text{-H%
\"{o}lder}}<\infty .$ In fact, we can find Lipschitz paths $x_{n}$ such that%
\begin{equation*}
\sup_{n}\left\Vert S_{N}\left( x_{n}\right) \right\Vert _{1/p\text{-H\"{o}%
lder}}\leq 3\left\Vert \mathbf{x}\right\Vert _{1/p\text{-H\"{o}lder}}.%
\newline
\end{equation*}%
(iii) Assume $p>1.$ Define $C^{0,1/p\text{-H\"{o}lder}}\left( \left[ 0,1%
\right] ,G^{N}\left( \mathbb{R}^{d}\right) \right) $ as the closure of
lifted Lipschitz paths $S_{N}\left( x\right) $ under the metric $d_{1/p\text{%
-H\"{o}lder}}$. For $\mathbf{x}\in C^{1/p\text{-H\"{o}lder}}\left( \left[ 0,1%
\right] ,G^{N}\left( \mathbb{R}^{d}\right) \right) $ we have%
\begin{equation*}
\mathbf{x}\in C^{0,1/p\text{-H\"{o}lder}}\left( \left[ 0,1\right]
,G^{N}\left( \mathbb{R}^{d}\right) \right) \text{ iff }r\left( \delta ;%
\mathbf{x}\right) \equiv \sup_{\substack{ 0\leq s<t\leq 1  \\ t-s\leq \delta 
}}\frac{\left\Vert \mathbf{x}_{s,t}\right\Vert }{\omega \left( s,t\right)
^{1/p}}\rightarrow 0\text{ as }\delta \rightarrow 0.
\end{equation*}%
In particular, $C^{0,1/p\text{-H\"{o}lder}}\left( \left[ 0,1\right]
,G^{N}\left( \mathbb{R}^{d}\right) \right) \subsetneq C^{1/p\text{-H\"{o}lder%
}}\left( \left[ 0,1\right] ,G^{N}\left( \mathbb{R}^{d}\right) \right) $.
\end{theorem}

One can see that $C^{0,1/p\text{-H\"{o}lder}}\left( \left[ 0,1\right]
,G^{N}\left( \mathbb{R}^{d}\right) \right) $ is Polish whereas $C^{1/p\text{%
-H\"{o}lder}}\left( \left[ 0,1\right] ,G^{N}\left( \mathbb{R}^{d}\right)
\right) $ lacks separability. Recall that $\left[ p\right] $ denotes the
integer part of some (positive) real number $p.$

\begin{definition}[\protect\cite{Ly, FV04}]
A path in $C^{1/p\text{-H\"{o}lder}}\left( \left[ 0,1\right] ,G^{\left[ p%
\right] }\left( \mathbb{R}^{d}\right) \right) $ is called a weak geometric $%
p $-rough path (with H\"{o}lder-control $\omega $). A path in $C^{0,1/p\text{%
-H\"{o}lder}}\left( \left[ 0,1\right] ,G^{\left[ p\right] }\left( \mathbb{R}%
^{d}\right) \right) $ is called a geometric $p$-rough path (with H\"{o}lder-
control $\omega $).
\end{definition}

\begin{proposition}[\protect\cite{Ly, LQ}]
\label{th1}Let $\mathbf{x}\in C^{1/p\text{-H\"{o}lder}}\left( \left[ 0,1%
\right] ,G^{\left[ p\right] }\left( \mathbb{R}^{d}\right) \right) $ and $N$ $%
>[p].$Then there exists a unique lift of $\mathbf{x}$ to a $G^{N}\left( 
\mathbb{R}^{d}\right) $-valued $1/p$-H\"{o}lder continuous path w.r.t.
Carnot-Caratheodory metric on $G^{N}\left( \mathbb{R}^{d}\right) $, denoted
by $S_{N}\left( \mathbf{x}\right) :\left[ 0,1\right] \rightarrow G^{N}\left( 
\mathbb{R}^{d}\right) $. Moreover, there exists a const $C_{\ref{th1}%
}=C\left( p,N\right) $ such that 
\begin{equation*}
\left\Vert S_{N}\left( \mathbf{x}\right) \right\Vert _{1/p\text{-H\"{o}lder;}%
\left[ 0,1\right] }\leq C_{\ref{th1}}\left\Vert \mathbf{x}\right\Vert _{1/p%
\text{-H\"{o}lder;}\left[ 0,1\right] }.
\end{equation*}
\end{proposition}

\section{Generalized Davie Estimates}

In this section we show that the step-$N$ Euler approximation is a good
approximation to ODE\ solutions in small time, even if we control only the
homogenous $1/p$-H\"{o}lder norm of $S_{N}\left( x\right) .$ In the case of $%
N=1,2$ this result is due to A. M. Davie, \cite{Da}. The existence of
geodesics associated to the Carnot-Caratheodory metric is our main tool to
generalize his results to the step-$N$ case.

\bigskip

Recall that a control ODE\ driven by $\mathrm{Lip}^{N}$ vector fields has
the step-$N$ Euler approximation%
\begin{equation*}
\pi \left( s,y_{s};x\right) _{s,t}\approx I\left[ y_{s},N,S_{N}\left(
x\right) _{s,t}\right] .
\end{equation*}

\bigskip The Geodesic Existence theorem, applied to $\mathbf{g=}$ $%
S_{N}\left( x\right) _{s,t}$, yields the shortest path in $\mathbb{R}^{d}$
whose iterated integrals mimick the first $N$ iterated integrals of the path
segement $x|_{[s,t]}.$ We called this path $x^{s,t}=\left(
x_{u}^{s,t}\right) _{u\in \left[ s,t\right] }.$ \ By construction, its step-$%
N$ Euler approximation over $\left[ s,t\right] $ is exactly equal to $%
I^{y_{s},N,S_{N}\left( x\right) _{s,t}}$ and we are led to the equally good 
\textit{step-}$N$\textit{\ geodesic approximation}%
\begin{equation*}
\pi \left( s,y_{s};x\right) _{s,t}\approx \pi \left( s,y_{s};x^{s,t}\right)
_{s,t}.
\end{equation*}%
This step-$N$ approximation is sometimes easier to handle. It also respect
the geometry given by the vector fields. Below, we shall use both. As last
preparation for the main result of this section, we need to understand the
regularity $y\mapsto I\left[ y,N,\mathbf{g}\right] $.

\begin{lemma}
\label{Eulerchangeinitialpoint}Assume that $\left( V_{i}\right) _{1\leq
i\leq d}\in \mathrm{Lip}^{N}\left( \mathbb{R}^{e}\right) $. For an element $%
\mathbf{g}\in G^{N}\left( \mathbb{R}^{d}\right) ,$ 
\begin{equation*}
\left\vert I^{y,N,\mathbf{g}}-I^{\tilde{y},N,\mathbf{g}}\right\vert \leq C_{%
\ref{Eulerchangeinitialpoint}}\left\vert y-\tilde{y}\right\vert \left(
\left\Vert \mathbf{g}\right\Vert +\left\Vert \mathbf{g}\right\Vert
^{N}\right)
\end{equation*}%
where $C_{\ref{Eulerchangeinitialpoint}}$ depends on $N$ and the $\mathrm{Lip%
}^{N}$ norm of the vector fields$.$
\end{lemma}

\begin{proof}
By definition of the Euler approximation $I^{y,N,\mathbf{g}}$%
\begin{equation*}
I^{y,N,\mathbf{g}}-I^{\tilde{y},N,\mathbf{g}}=\sum_{k=1}^{N}\sum_{\substack{ %
i_{1},\cdots ,i_{k}  \\ \in \left\{ 1,\cdots ,d\right\} }}[V_{i_{1}}\cdots
V_{i_{k}}H\left( y\right) -V_{i_{1}}\cdots V_{i_{k}}H\left( \tilde{y}\right)
]\mathbf{g}^{k,i_{1},\cdots ,i_{k}}.
\end{equation*}%
Since $y\mapsto V_{i_{1}}\cdots V_{i_{k}}H\left( y\right) $ is Lipschitz,%
\begin{equation*}
\left\vert I^{y,N,\mathbf{g}}-I^{\tilde{y},N,\mathbf{g}}\right\vert \leq C_{%
\ref{Eulerchangeinitialpoint}}^{1}\sum_{k=1}^{N}\sum_{\substack{ %
i_{1},\cdots ,i_{k}  \\ \in \left\{ 1,\cdots ,d\right\} }}\left\vert y-%
\tilde{y}\right\vert \left\vert \mathbf{g}^{k,i_{1},\cdots
,i_{k}}\right\vert .
\end{equation*}%
From equivalence of homogenous norms, $\left\vert \mathbf{g}^{k,i_{1},\cdots
,i_{k}}\right\vert \leq C_{\ref{Eulerchangeinitialpoint}}^{2}\left\Vert 
\mathbf{g}\right\Vert ^{k}$ and hence 
\begin{equation*}
\left\vert I^{y,N,\mathbf{g}}-I^{\tilde{y},N,\mathbf{g}}\right\vert \leq C_{%
\ref{Eulerchangeinitialpoint}}^{3}\left\vert y-\tilde{y}\right\vert \max
\left\{ \left\Vert \mathbf{g}\right\Vert ^{N},\left\Vert \mathbf{g}%
\right\Vert \right\} .
\end{equation*}
\end{proof}

The next lemma is technical but very important. It quantifies the quality of
step-$N$ Euler and geodesic approximations and gives ODE\ bounds which do
not blow up with the Lipschitz norm of the driving signal. Recall that $%
\omega \left( s,t\right) \equiv t-s$ although the proof can adapted to
general super additive control function \cite{LQ}.

\begin{lemma}[Generalized Davie Lemma]
\label{tfa}Let $p\geq 1\ $and $\left( V_{i}\right) _{1\leq i\leq d}\in 
\mathrm{Lip}^{N}\left( \mathbb{R}^{e}\right) $ for some integer $N>p-1$.
Assume that\newline
(i) $x:\left[ 0,1\right] \rightarrow \mathbb{R}^{d}$ is a Lipschitz path
with step-$N$ lift $\mathbf{x}=S_{N}\left( x\right) $ and $\left\Vert 
\mathbf{x}\right\Vert _{1/p\text{-H\"{o}lder}}\leq M_{1}$.\newline
(ii) $y_{0}\in \mathbb{R}^{e}$ with $\left\vert y_{0}\right\vert \leq M_{2};$
Then, there exists a positive constant $C_{\ref{tfa}}=C_{\ref{tfa}}\left(
M_{1}\right) $, also dependent on $p,N,M_{2}$ and the vector fields $%
V_{1},...,V_{d}$ but \textbf{not} dependent on the Lipschitz norm of \ $x$,
such that for all $0\leq s\leq t\leq 1$%
\begin{equation}
\left\vert \pi \left( 0,y_{0};x\right) _{s,t}\right\vert <C_{\ref{tfa}%
}\left( M_{1}\right) \omega \left( s,t\right) ^{\frac{1}{p}},
\label{Holdernormofsolution}
\end{equation}%
and with $\theta =\frac{N+1}{p}>1,$%
\begin{equation}
\left\vert \pi \left( 0,y_{0};x\right) _{s,t}-I^{y_{s},N,\mathbf{x}%
_{s,t}}\right\vert \leq C_{\ref{tfa}}\left( M_{1}\right) \omega \left(
s,t\right) ^{\theta }.  \label{ApproximationbyEuler}
\end{equation}%
Moreover, if we assume\newline
(iii) $x^{s,t}:\left[ s,t\right] \rightarrow \mathbb{R}^{d}$ are Lipschitz
paths such that $S_{N}\left( x^{s,t}\right) _{s,t}=\mathbf{x}_{s,t}$ and
such that%
\begin{equation}
\int_{s}^{t}\left\vert dx^{s,t}\right\vert \leq KM_{1}\omega \left(
s,t\right) ^{1/p}  \label{XST}
\end{equation}%
for some positive real $K$ then%
\begin{equation}
\left\vert \pi \left( 0,y_{0};x\right) _{s,t}-\pi \left(
s,y_{s};x^{s,t}\right) _{s,t}\right\vert \leq C_{\ref{tfa}}^{\prime }\left(
M_{1}\right) \omega \left( s,t\right) ^{\theta },
\label{Approximationbygeodesics}
\end{equation}%
with $C_{\ref{tfa}}^{\prime }=C_{\ref{tfa}}^{\prime }\left( M_{1}\right) $
also dependent on $p,N,M_{2},K$ and the vector fields $V_{1},...,V_{d}$.
\end{lemma}

\begin{remark}
By hypothesis, $\left\Vert \mathbf{x}_{s,t}\right\Vert \leq M_{1}\omega
\left( s,t\right) ^{1/p}$ and we can find a (projected) geodesic $\gamma
^{\ast }$:$\left[ s,t\right] \rightarrow \mathbb{R}^{d}$ with signature $%
\mathbf{x}_{s,t}$ and length $\int_{s}^{t}\left\vert d\gamma ^{\ast
}\right\vert =\left\Vert \mathbf{x}_{s,t}\right\Vert \leq $ $M_{1}\omega
\left( s,t\right) ^{1/p}$. In other words, paths $\left\{ x^{s,t}\right\} $
as postulated in (iii) always exist, even for $K=1$. (By not fixing $K$ we
get a more general approximation result and highlight what is needed in the
proof.)
\end{remark}

\begin{proof}
Without loss of generality, we assume that $M_{1}\geq 1$ (otherwise set $%
M_{1}=1$.)$.$

Write $y_{t}=\pi \left( 0,y_{0};x\right) _{t}$ and $\Gamma
_{s,t}=y_{s,t}-\pi \left( s,y_{s};x^{s,t}\right) _{s,t}.$ We first show (\ref%
{Approximationbygeodesics}) and divide the argument in two steps.

\underline{\textit{First Step:}} Fix $0\leq s<t<u\leq 1$. We try to control $%
\Gamma _{s,u}$ in terms of $\Gamma _{s,t}$ and $\Gamma _{t,u}$. To this end,
it is useful to define $x^{s,t,u}$ to be the concatenation of $x^{s,t}$ and $%
x^{t,u}.$ Observe that $\left. x^{s,t,u}\right\vert _{\left[ s,u\right] }$
has the step-$N$ signature $\mathbf{x}_{s,t}\otimes \mathbf{x}_{t,u}=\mathbf{%
x}_{s,u}$ and 
\begin{equation}
\int_{s}^{u}\left\vert dx^{s,t,u}\right\vert =\int_{s}^{t}\left\vert
dx^{s,t}\right\vert +\int_{t}^{u}\left\vert dx^{t,u}\right\vert \leq
2KM_{1}\omega \left( s,u\right) ^{1/p}.  \label{XSTU}
\end{equation}%
By uniqueness of ODE\ solutions,%
\begin{equation*}
\pi \left( s,y_{s};x^{s,t,u}\right) _{r}=\left\{ 
\begin{array}{l}
\pi \left( s,y_{s};x^{s,t}\right) _{r}\text{ if }r\in \left[ s,t\right] \\ 
\pi \left( t,\pi \left( s,y_{s};x^{s,t}\right) _{r};x^{t,u}\right) _{r}\text{
if }r\in \left[ t,u\right] .%
\end{array}%
\right.
\end{equation*}%
We have 
\begin{eqnarray*}
-\Gamma _{s,u}+\Gamma _{s,t}+\Gamma _{t,u} &=&\pi \left(
s,y_{s};x^{s,u}\right) _{s,u}-\pi \left( s,y_{s};x^{s,t}\right) _{s,t}-\pi
\left( t,y_{t};x^{t,u}\right) _{t,u} \\
&=&\pi \left( s,y_{s};x^{s,u}\right) _{s,u}-\pi \left(
s,y_{s};x^{s,t,u}\right) _{s,u} \\
&&+\pi \left( s,y_{s};x^{s,t,u}\right) _{s,u}-\pi \left(
s,y_{s};x^{s,t}\right) _{s,t}-\pi \left( t,y_{t};x^{t,u}\right) _{t,u}
\end{eqnarray*}%
\newline
By defintion of $\pi \left( s,y_{s};x^{s,t,u}\right) $,%
\begin{equation*}
\pi \left( s,y_{s};x^{s,t,u}\right) _{s,u}-\pi \left( s,y_{s};x^{s,t}\right)
_{s,t}=\pi \left( t,\pi \left( s,y_{s};x^{s,t}\right) _{t};x^{t,u}\right)
_{t,u},
\end{equation*}%
hence,%
\begin{eqnarray*}
-\Gamma _{s,u}+\Gamma _{s,t}+\Gamma _{t,u} &=&\pi \left(
s,y_{s};x^{s,u}\right) _{s,u}-\pi \left( s,y_{s};x^{s,t,u}\right) _{s,u} \\
&&+\pi \left( t,\pi \left( s,y_{s};x^{s,t}\right) _{t};x^{t,u}\right)
_{t,u}-\pi \left( t,y_{t};x^{t,u}\right) _{t,u}.
\end{eqnarray*}%
In particular, from lemma \ref{Euler}, we have%
\begin{multline*}
\left\vert \pi \left( s,y_{s};x^{s,u}\right) _{s,u}-\pi \left(
s,y_{s};x^{s,t,u}\right) _{s,u}\right\vert \\
\left. 
\begin{array}{l}
\leq \left\vert \pi \left( s,y_{s};x^{s,u}\right) _{s,u}-I^{y_{s},N,\mathbf{x%
}_{s,u}}\right\vert +\left\vert \pi \left( s,y_{s};x^{s,t,u}\right)
_{s,u}-I^{y_{s},N,\mathbf{x}_{s,u}}\right\vert \\ 
\leq 2C_{\ref{tfa}}^{1}\left( M_{1}\right) \omega \left( s,u\right) ^{\frac{%
N+1}{p}}=2C_{\ref{tfa}}^{1}\left( M_{1}\right) \omega \left( s,u\right)
^{\theta }.%
\end{array}%
\right.
\end{multline*}%
using (\ref{XST}) and (\ref{XSTU}). Then,%
\begin{multline*}
\left\vert \pi \left( t,\pi \left( s,y_{s};x^{s,t}\right)
_{t};x^{t,u}\right) _{t,u}-\pi \left( t,y_{t};x^{t,u}\right)
_{t,u}\right\vert \\
\left. 
\begin{array}{l}
\leq \left\vert \pi \left( t,\pi \left( s,y_{s};x^{s,t}\right)
_{t};x^{t,u}\right) _{t,u}-I^{\pi \left( s,y_{s};x^{s,t}\right) _{t},N,%
\mathbf{x}_{t,u}}\right\vert +\left\vert \pi \left( t,y_{t};x^{t,u}\right)
_{t,u}-I^{y_{t},N,\mathbf{x}_{t,u}}\right\vert \\ 
\text{ \ \ }+\left\vert I^{\pi \left( s,y_{s};x^{s,t}\right) _{t},N,\mathbf{x%
}_{t,u}}-I^{y_{t},N,\mathbf{x}_{t,u}}\right\vert .%
\end{array}%
\right.
\end{multline*}%
Once again, by lemma \ref{Euler}, 
\begin{multline*}
\left\vert \pi \left( t,\pi \left( s,y_{s};x^{s,t}\right)
_{t};x^{t,u}\right) _{t,u}-I^{\pi \left( s,y_{s};x^{s,t}\right) _{t},N,%
\mathbf{x}_{t,u}}\right\vert +\left\vert \pi \left( t,y_{t};x^{t,u}\right)
_{t,u}-I^{y_{t},N,\mathbf{x}_{t,u}}\right\vert \\
\leq 2C_{\ref{tfa}}^{2}\left( M_{1}\right) \omega \left( t,u\right) ^{\frac{%
N+1}{p}}\leq 2C_{\ref{tfa}}^{2}\left( M_{1}\right) \omega \left( s,u\right)
^{\theta }.
\end{multline*}%
Finally, by lemma \ref{Eulerchangeinitialpoint},%
\begin{eqnarray*}
\left\vert I^{\pi \left( s,y_{s};x^{s,t}\right) _{t},N,\mathbf{x}%
_{t,u}}-I^{y_{t},N,\mathbf{x}_{t,u}}\right\vert &\leq &C_{\ref{tfa}%
}^{3}\left\vert \Gamma _{s,t}\right\vert \max \left( M_{1}\omega \left(
t,u\right) ^{1/p},M_{1}^{N}\omega \left( t,u\right) ^{N/p}\right) \\
&\leq &C_{\ref{tfa}}^{3}\left\vert \Gamma _{s,t}\right\vert M_{1}^{N}\max
\left( \omega \left( t,u\right) ^{1/p},\omega \left( t,u\right) ^{N/p}\right)
\\
&\leq &2C_{\ref{tfa}}^{3}\left( M_{1}\right) \left\vert \Gamma
_{s,t}\right\vert \omega \left( t,u\right) ^{1/p}
\end{eqnarray*}%
using that $\omega \left( t,u\right) \leq 1$. Putting the pieces together,
we have%
\begin{equation*}
\left\vert -\Gamma _{s,u}+\Gamma _{s,t}+\Gamma _{t,u}\right\vert \leq \left(
2C_{\ref{tfa}}^{1}\left( M_{1}\right) +2C_{\ref{tfa}}^{2}\left( M_{1}\right)
\right) \omega \left( s,u\right) ^{\theta }+2C_{\ref{tfa}}^{3}\left(
M_{1}\right) \left\vert \Gamma _{s,t}\right\vert \omega \left( t,u\right)
^{1/p}.
\end{equation*}%
It follows that%
\begin{eqnarray}
\left\vert \Gamma _{s,u}\right\vert &\leq &\left\vert -\Gamma _{s,u}+\Gamma
_{s,t}+\Gamma _{t,u}\right\vert +\left\vert \Gamma _{s,t}\right\vert
+\left\vert \Gamma _{t,u}\right\vert  \notag \\
&\leq &\left\vert \Gamma _{s,t}\right\vert \left( 1+C_{\ref{tfa}}^{4}\left(
M_{1}\right) \omega \left( t,u\right) ^{1/p}\right) +\left\vert \Gamma
_{t,u}\right\vert +C_{\ref{tfa}}^{4}\left( M_{1}\right) \omega \left(
s,u\right) ^{\theta }  \label{inGamma}
\end{eqnarray}%
\underline{\textit{Second Step:}} For $0\leq s<t<u\leq 1$ inequality (\ref%
{inGamma}) can be rewritten as%
\begin{eqnarray*}
\frac{\left\vert \Gamma _{s,u}\right\vert }{\omega \left( s,u\right)
^{\theta }} &\leq &\frac{\omega \left( s,t\right) ^{\theta }}{\omega \left(
s,u\right) ^{\theta }}\frac{\left\vert \Gamma _{s,t}\right\vert }{\omega
\left( s,t\right) ^{\theta }}\left( 1+C_{\ref{tfa}}^{4}\left( M_{1}\right)
\omega \left( t,u\right) ^{1/p}\right) \\
&&+\frac{\omega \left( t,u\right) ^{\theta }}{\omega \left( s,u\right)
^{\theta }}\frac{\left\vert \Gamma _{t,u}\right\vert }{\omega \left(
t,u\right) ^{\theta }}+C_{\ref{tfa}}^{4}\left( M_{1}\right) .
\end{eqnarray*}%
Define for $r\in (0,1]$,%
\begin{equation*}
\varrho \left( r\right) =\sup_{\substack{ 0\leq s<t\leq 1  \\ \omega \left(
s,t\right) \leq r}}\frac{\left\vert \Gamma _{s,t}\right\vert }{\omega \left(
s,t\right) ^{\theta }}.
\end{equation*}%
Note that $\rho \left( r\right) <\infty $. Indeed, this follows from%
\begin{eqnarray*}
\left\vert \Gamma _{s,t}\right\vert &=&\left\vert y_{s,t}-\pi \left(
s,y_{s};x^{s,t}\right) _{s,t}\right\vert \\
&\leq &\left\vert y_{s,t}-I^{y_{s},N,S_{N}\left( x\right) _{s,t}}\right\vert
+\left\vert \pi \left( s,y_{s};x^{s,t}\right) _{s,t}-I^{y_{s},N,S_{N}\left(
x\right) _{s,t}}\right\vert \\
&\leq &C_{\ref{tfa}}^{5}\left( \int_{s}^{t}\left\vert dx_{r}\right\vert
\right) ^{N+1}+C_{\ref{tfa}}^{5}\left( \int_{s}^{t}\left\vert
dx_{r}^{s,t}\right\vert \right) ^{N+1} \\
&\leq &C_{\ref{tfa}}^{5}\left\vert x\right\vert _{Lip}^{N+1}\left\vert
t-s\right\vert ^{N+1}+C_{\ref{tfa}}^{5}KM_{1}\omega \left( s,t\right) ^{%
\frac{N+1}{p}}.
\end{eqnarray*}%
The problem with this bound is that it blows up with $\left\vert
x\right\vert _{Lip}$. The argument which follows shows that, in fact, $\rho
\left( r\right) $ will \textit{not} blow up with $\left\vert x\right\vert
_{Lip}.$ Pick arbitrary points $s<u$ such that $\omega \left( s,u\right)
\leq r$, amd set $t=\left( s+u\right) /2$ so that%
\begin{equation*}
\omega \left( s,t\right) =\omega \left( t,u\right) =\frac{1}{2}\omega \left(
s,u\right) .
\end{equation*}%
We obtain from inequality (\ref{inGamma}) that%
\begin{eqnarray*}
\frac{\left\vert \Gamma _{s,u}\right\vert }{\omega \left( s,u\right)
^{\theta }} &\leq &\left( \frac{1}{2}\right) ^{\theta }\varrho \left(
r/2\right) \left( 1+C_{\ref{tfa}}^{4}\left( M_{1}\right) r^{1/p}\right) \\
&&+\left( \frac{1}{2}\right) ^{\theta }\varrho \left( r/2\right) +C_{\ref%
{tfa}}^{4}\left( M_{1}\right) \\
&\leq &2^{1-\theta }\varrho \left( r/2\right) \underset{\equiv \gamma \left(
r\right) }{\underbrace{\left( 1+C_{\ref{tfa}}^{4}\left( M_{1}\right)
r^{1/p}\right) }}+C_{\ref{tfa}}^{4}\left( M_{1}\right)
\end{eqnarray*}%
and taking the supremum over all $s<u$ with $\omega \left( s,u\right) \leq r$
gives%
\begin{equation*}
\varrho \left( r\right) \leq 2^{1-\theta }\varrho \left( r/2\right) \gamma
\left( r\right) +C_{\ref{tfa}}^{4}\text{.}
\end{equation*}%
After $n$ iterated uses of the inequality for $\rho $ we find%
\begin{eqnarray*}
\varrho \left( r\right) &\leq &\left( 2^{1-\theta }\right)
^{n+1}\prod_{k=0}^{n}\gamma \left( \frac{r}{2^{k}}\right) \varrho \left( 
\frac{r}{2^{n+1}}\right) \\
&&+C_{\ref{tfa}}^{4}\left( M_{1}\right) \left[ \sum_{k=0}^{n}\left( \left(
2^{1-\frac{N+1}{p}}\right) ^{k}\prod_{j=0}^{k-1}\gamma \left( \frac{r}{2^{j}}%
\right) \right) \right] .
\end{eqnarray*}%
Let $C_{\ref{tfa}}^{5}\left( r,M_{1}\right) :=\sum_{k=0}^{\infty }\left(
\left( 2^{1-\frac{N+1}{p}}\right) ^{k}\prod_{j=0}^{k-1}\gamma \left( \frac{r%
}{2^{j}}\right) \right) .$ Note that $\prod_{k=0}^{n}\gamma \left(
r/2^{k}\right) $ is increasing in $n$ and since $\gamma \left( r\right) \leq
e^{C_{\ref{tfa}}^{4}\left( M_{1}\right) r^{1/p}},$ the supremum over $r$ of
the infinite product $\prod_{k=0}^{\infty }\gamma \left( r/2^{k}\right) $ is
finite, which implies that $\sup_{0\leq r\leq 1}C_{\ref{tfa}}^{5}\left(
r,M_{1}\right) $ is also finite.

Hence,%
\begin{equation*}
\varrho \left( r\right) \leq C_{\ref{tfa}}^{5}\left( M\right) \left(
2^{1-\theta }\right) ^{n+1}\varrho \left( \frac{r}{2^{n+1}}\right) +\frac{C_{%
\ref{tfa}}^{4}\left( M_{1}\right) C_{\ref{tfa}}^{5}\left( r,M_{1}\right) }{%
1-2^{1-\theta }}
\end{equation*}%
and sending $n\rightarrow \infty $ leaves us with (note $\theta >1$ here),%
\begin{equation*}
\varrho \left( r\right) \leq \frac{C_{\ref{tfa}}^{4}\left( M_{1}\right) C_{%
\ref{tfa}}^{5}\left( r,M_{1}\right) }{1-2^{1-\theta }}.
\end{equation*}%
From the very definition of $\rho $ with $r=1$ we obtain%
\begin{equation*}
\left\vert \Gamma _{s,t}\right\vert \leq C_{\ref{tfa}}^{6}\left(
M_{1}\right) \omega \left( s,t\right) ^{\theta }.
\end{equation*}

\underline{\textit{Third Step:}} Using lemma \ref{boundY} and (\ref{XST}),%
\begin{eqnarray*}
\left\vert \pi \left( s,y_{s};x^{s,t}\right) _{s,t}\right\vert &\leq &C_{\ref%
{tfa}}^{7}\int_{s}^{t}\left\vert dx^{s,t}\right\vert \\
&\leq &C_{\ref{tfa}}^{7}KM_{1}\omega \left( s,t\right) ^{1/p}.
\end{eqnarray*}%
Then, for all $s,t\in \left[ 0,1\right] $%
\begin{eqnarray*}
\left\vert y_{s,t}\right\vert &\leq &\left\vert y_{s,t}-\pi \left(
s,y_{s};x^{s,t}\right) _{s,t}\right\vert +\left\vert \pi \left(
s,y_{s};x^{s,t}\right) _{s,t}\right\vert \\
&\leq &C_{\ref{tfa}}^{6}\left( M_{1}\right) \omega \left( s,t\right)
^{\theta }+C_{\ref{tfa}}^{7}KM_{1}\omega \left( s,t\right) ^{1/p} \\
&\leq &C_{\ref{tfa}}^{6}\left( M_{1}\right) \omega \left( s,t\right)
^{1/p}+C_{\ref{tfa}}^{7}KM_{1}\omega \left( s,t\right) ^{1/p} \\
&\equiv &C_{\ref{tfa}}^{8}\left( M_{1}\right) \omega \left( s,t\right)
^{1/p}.
\end{eqnarray*}%
Although $C_{\ref{tfa}}^{8}$ manifestly depends on $K$, we may specialize
the construction using geodesics $\left\{ x^{s,t}\right\} $for which $K=1$.
With such paths, $C_{\ref{tfa}}^{8}\left( M_{1}\right) $ of course would not
depend on $K.$ In particular, the H\"{o}lder norm on $y$ does not depend on $%
K.$\newline
\underline{\textit{Fourth Step:}} Finally, (\ref{ApproximationbyEuler}) is
obtained from \ (\ref{Approximationbygeodesics}) via triangle inequality and
lemma \ref{Euler}, taking into account (\ref{XST}).
\end{proof}

\begin{corollary}
\label{corTFA}There exists a constant $C_{\ref{corTFA}},$ which may depend
on $p,N,M_{2}$ and the vector fields $V_{1},...,V_{d}$ so that for all $%
M_{1}\geq 1$ 
\begin{equation*}
C_{\ref{tfa}}\leq C_{\ref{corTFA}}\exp \left( 12N^{2}\ln \left( M_{1}\right)
^{2}\right) ,
\end{equation*}%
This implies (the $O$-notation being understood as $M_{1}\rightarrow \infty $%
)%
\begin{equation*}
\ln C_{\ref{tfa}}\left( M_{1}\right) =O\left( \left( \ln M_{1}\right)
^{2}\right) .
\end{equation*}%
The same estimates holds for $C_{\ref{tfa}}^{\prime }$, allowing for
additional dependence on $K$.
\end{corollary}

\begin{proof}
Inspection of the first step in the proof of Davie's lemma shows that 
\begin{equation*}
C_{\ref{tfa}}^{1}\left( M_{1}\right) ,C_{\ref{tfa}}^{2}\left( M_{1}\right)
,C_{\ref{tfa}}^{3}\left( M_{1}\right) ,C_{\ref{tfa}}^{4}\left( M_{1}\right)
\leq C_{\ref{corTFA}}^{1}M_{1}^{N+1}\text{.}
\end{equation*}%
The only difficulty is to control%
\begin{eqnarray*}
C_{\ref{tfa}}^{5}\left( 1,M_{1}\right) &=&1+\sum_{k=1}^{\infty }\left(
\left( 2^{1-\frac{N+1}{p}}\right) ^{k}\prod_{j=0}^{k-1}\gamma \left( \frac{1%
}{2^{j}}\right) \right) \\
&\leq &\sum_{k=0}^{\infty }\left( \left( 2^{1-\frac{N+1}{p}}\right)
^{k}\prod_{j=0}^{k-1}\left( 1+C_{\ref{corTFA}}^{1}M_{1}^{N+1}2^{-jp}\right)
\right) .
\end{eqnarray*}%
To understand the dependence of the right hand side on $M_{1},$ we define
the function for some fixed $a\in \left( 0,1\right) $,%
\begin{eqnarray*}
\Lambda \left( k,b\right) &=&a^{k}\Pi _{j=0}^{k-1}\left( 1+b2^{-pj}\right) \\
\Gamma \left( n,b\right) &=&\sum_{k=0}^{n}\Lambda \left( k,b\right) .\text{ }
\end{eqnarray*}%
We need to understand the dependence of $\lim_{n\rightarrow \infty }\Gamma
\left( n,b\right) $ on $b$. One could use a naive approach (the one used in
the previous proof) to get%
\begin{eqnarray}
\Lambda \left( k,b\right) &\leq &a^{k}\Pi _{j=0}^{\infty -1}\left(
1+b2^{-pj}\right)  \label{simpleEstimate} \\
&\leq &a^{k}\Pi _{j=0}^{\infty -1}\exp \left( b2^{-pj}\right)  \notag \\
&\leq &a^{k}\exp \left( \frac{b}{1-2^{-p}}\right) ,  \notag
\end{eqnarray}%
and hence, 
\begin{equation*}
\lim_{n\rightarrow \infty }\Gamma \left( n,b\right) \leq \frac{1}{1-a}\exp
\left( \frac{b}{1-2^{-p}}\right) .
\end{equation*}%
Unfortunately, the right hand side in the last equation grows to fast in $b$
for our purposes. To obtain a better estimate, we first observe that 
\begin{equation*}
\Lambda \left( k,0\right) =a^{k}\text{ and }\frac{\partial }{\partial b}%
\Lambda \left( k,b\right) =\left( \sum_{j=0}^{k-1}\frac{2^{-pj}}{1+b2^{-pj}}%
\right) \Lambda \left( k,b\right) .
\end{equation*}%
Then we note that $x\mapsto $ $2^{-px}/\left( 1+b2^{-px}\right)
=b^{-1}\left( 1-\left( 1+b2^{-px}\right) ^{-1}\right) $ is decreasing in $x$
so that%
\begin{eqnarray*}
\sum_{j=0}^{k-1}\frac{2^{-pj}}{1+b2^{-pj}} &\leq &\frac{1}{1+b}%
+\sum_{j=1}^{\infty }\frac{2^{-pj}}{1+b2^{-pj}} \\
&\leq &\frac{1}{1+b}+\int_{0}^{\infty }\frac{2^{-px}dx}{1+b2^{-px}} \\
&=&\frac{1}{1+b}+\frac{\ln \left( 1+b\right) }{b\ln \left( 2^{p}\right) }%
\text{ \ \ } \\
&\leq &3\frac{\ln b}{b}
\end{eqnarray*}%
for all~$b\geq e$. We also note that 
\begin{equation*}
f\left( b\right) :=\exp \left( \frac{3}{2}\left( \ln b\right) ^{2}\right) 
\text{ solves }\frac{\partial }{\partial b}f\left( b\right) =3\frac{\ln b}{b}%
f\left( b\right) \text{. }
\end{equation*}%
and any other solution to this ODE must be a multiple of $f$. By ODE
comparison we see that, for $b\geq e$,%
\begin{equation*}
\frac{\Lambda \left( k,b\right) }{\Lambda \left( k,e\right) }\leq \frac{%
f\left( b\right) }{f\left( e\right) }\leq e^{\frac{3}{2}\left( \ln b\right)
^{2}},
\end{equation*}%
which implies that, using (\ref{simpleEstimate}) 
\begin{eqnarray*}
\Lambda \left( k,b\right) &\leq &\Lambda \left( k,e\right) e^{\frac{3}{2}%
\left( \ln b\right) ^{2}} \\
&\leq &a^{k}\exp \left( \frac{e}{1-2^{-p}}\right) e^{\frac{3}{2}\left( \ln
b\right) ^{2}}
\end{eqnarray*}%
After summing over all non-negative integers $k$ we see that%
\begin{equation*}
\lim_{n\rightarrow \infty }\Gamma \left( n,b\right) \leq \frac{\exp \left( 
\frac{e}{1-2^{-p}}\right) }{1-a}e^{\frac{3}{2}\left( \ln b\right) ^{2}}.
\end{equation*}%
Hence, we have proved that 
\begin{eqnarray*}
C_{\ref{tfa}}^{5}\left( 1,M_{1}\right) &\leq &\frac{\exp \left( \frac{e}{%
1-2^{-p}}\right) }{1-2^{1-\frac{N+1}{p}}}\exp \left( \frac{3}{2}\ln \left(
C_{\ref{corTFA}}^{1}M_{1}^{N+1}\right) ^{2}\right) \\
&\leq &C_{\ref{corTFA}}^{2}\exp \left( 3\left( N+1\right) ^{2}\ln \left(
M_{1}\right) ^{2}\right) \\
&\leq &C_{\ref{corTFA}}^{2}\exp \left( 6N^{2}\ln \left( M_{1}\right)
^{2}\right) .
\end{eqnarray*}%
This lead to 
\begin{eqnarray*}
C_{\ref{tfa}}\left( M_{1}\right) &\leq &C_{\ref{corTFA}}^{3}M_{1}^{N+1}\exp
\left( 6N^{2}\ln \left( M_{1}\right) ^{2}\right) \\
&\leq &C_{\ref{corTFA}}^{3}\exp \left( 12N^{2}\ln \left( M_{1}\right)
^{2}\right) ,\text{ for }M_{1}\geq 3.
\end{eqnarray*}%
By increasing $C_{\ref{corTFA}}^{3}$ if needed we can assume that this
estimate holds for all $M_{1}\geq 1$. Clearly, the same estimate holds for $%
C_{\ref{tfa}}^{\prime }\left( M_{1}\right) .$
\end{proof}

\section{Euler Estimates for Rough Differential Equations (RDEs)}

We consider controlled differential equations in the sense of T. Lyons. The
driving signal is assumed to be a weak geometric $p$-rough path with H\"{o}%
lder control $\omega \left( s,t\right) =t-s$. Recall that this means $%
\mathbf{x}:\left[ 0,1\right] \rightarrow G^{\left[ p\right] }\left( \mathbb{R%
}^{d}\right) $ is $1/p$-H\"{o}lder continuous w.r.t. Carnot-Caratheodory
metric on $G^{\left[ p\right] }\left( \mathbb{R}^{d}\right) $. Lyons' theory 
\cite{Ly, LQ, Ly04} then implies existence and uniquess of a solution to the
differential equations driven by $\mathbf{x}$ along vector fields $%
V_{1},...\,V_{d}\in Lip^{p+\epsilon \text{ }}\left( \mathbb{R}^{e}\right) $
started at some point $y_{0}\in \mathbb{R}^{e}$ at time $0$. This \textit{%
RDE solution} is also a (weak) geometric $p$-rough path, over $\mathbb{R}%
^{e} $ instead of $\mathbb{R}^{d}$, denoted by 
\begin{equation*}
\mathbf{\pi }\left( 0,y_{0},\mathbf{x}\right) \equiv \mathbf{y}\text{,}
\end{equation*}%
with the same modulus of continuity as $\mathbf{x}$. For our application it
will be sufficient to consider the \textit{pathlevel RDE solution} (obtained
by projection)%
\begin{equation*}
\pi \left( 0,y_{0},\mathbf{x}\right) \equiv y:\left[ 0,1\right] \rightarrow 
\mathbb{R}^{e}.
\end{equation*}%
Thus, $y$ is an$\,\mathbb{R}^{e}$-valued $1/p$-H\"{o}lder continuous path in
the usual sense.

\begin{theorem}
\label{EulerRDE} Fix an integer $N\geq \lbrack p]+1$ and $\mathrm{Lip}^{N}$%
-vector fields $V_{1},...,V_{d}$ on $\mathbb{R}^{e}$. Let $\mathbf{x}$ be a
weak geometric $p$-rough path with $\left\Vert \mathbf{x}\right\Vert
_{p,\omega }\leq M$. Then there exists a unique pathlevel RDE solution $\pi
\left( 0,y_{0},\mathbf{x}\right) \equiv y$. Moreover, (a) there exists
constant $C_{\ref{EulerRDE}}=C_{\ref{EulerRDE}}\left( M\right) $, also
dependent on $p,N,y_{0}$ and $V_{1},...,V_{d}$, such that%
\begin{equation*}
\,\left\vert y\right\vert _{1/p\text{-H\"{o}lder;}\left[ 0,1\right] }\leq C_{%
\ref{EulerRDE}}
\end{equation*}%
and (b) a constant $C_{\ref{EulerRDE}}^{\prime }=C_{\ref{EulerRDE}}^{\prime
}\left( M\right) $ with similar dependencies such that for all $0\leq
s<t\leq 1$,%
\begin{equation*}
\left\vert y_{s,t}-I\left[ y_{s},N,S_{N}\left( \mathbf{x}\right) _{s,t}%
\right] \right\vert \leq C_{\ref{EulerRDE}}^{\prime }\omega \left(
s,t\right) ^{\theta }\text{ with }\theta =\frac{N+1}{p}>1\text{.}
\end{equation*}%
Finally, keeping all parameters but $M$ fixed,%
\begin{equation*}
\ln C_{\ref{EulerRDE}},\ln C_{\ref{EulerRDE}}^{\prime }=O\left( \left( \ln
M\right) ^{2}\right) \text{ as }M\rightarrow \infty \text{.}
\end{equation*}
\end{theorem}

\begin{proof}
$\mathrm{Lip}^{[p]+1}$-regularity is more than enough to ensure existence
and uniqueness of RDE solutions, see \cite{Ly, LQ, Ly04}.\newline
(a) From Theorem \ref{ThmOnHoelderRoughPathsProperties} we can find
Lipschitz paths $x^{n}$ such that%
\begin{equation*}
S_{\left[ p\right] }\left( x^{n}\right) \rightarrow \mathbf{x}
\end{equation*}%
uniformly on $\left[ 0,1\right] $, such that%
\begin{equation*}
\left\Vert S_{\left[ p\right] }\left( x^{n}\right) \right\Vert _{1/p\text{-H%
\"{o}lder}}\leq 3M\equiv M_{1}
\end{equation*}%
The Universal Limit Theorem implies a forteriori that 
\begin{equation*}
\pi \left( 0,y_{0};x^{n}\right) \rightarrow \pi \left( 0,y_{0};\mathbf{x}%
\right)
\end{equation*}%
uniformly on $\left[ 0,1\right] $. On the other hand, Davie's lemma implies
that%
\begin{equation*}
\sup_{n}\left\vert \pi \left( 0,y_{0};x^{n}\right) \right\vert _{1/p\text{-H%
\"{o}lder}}\leq C_{\ref{EulerRDE}}^{1}<\infty
\end{equation*}%
where $C_{\ref{EulerRDE}}^{1}$ is the constant $C_{\ref{tfa}}=C_{\ref{tfa}%
}\left( M_{1}\right) $ from lemma \ref{tfa}. It follows that%
\begin{equation*}
\left\vert \pi \left( 0,y_{0};\mathbf{x}\right) \right\vert _{1/p\text{-H%
\"{o}lder}}\leq C_{\ref{EulerRDE}}^{1}<\infty .
\end{equation*}%
From corollary \ref{corTFA},%
\begin{equation*}
\ln C_{\ref{EulerRDE}}^{1}=O\left( \left( \ln M_{1}\right) ^{2}\right)
=O\left( \left( \ln M\right) ^{2}\right) \text{.}
\end{equation*}

(b)\ By lemma \ref{th1}, a weak geometric $p$-rough path $\mathbf{x}$ with $%
\left\Vert \mathbf{x}\right\Vert _{1/p\text{-H\"{o}lder}}\leq M$ lifts
uniquely to a path $S_{N}\left( \mathbf{x}\right) \equiv \mathbf{\bar{x}}:%
\left[ 0,1\right] \rightarrow G^{N}\left( \mathbb{R}^{d}\right) $ such that%
\begin{equation*}
\left\Vert \mathbf{\bar{x}}\right\Vert _{1/p\text{-H\"{o}lder}}\leq C_{\ref%
{EulerRDE}}^{2}M
\end{equation*}%
for some constant $C_{\ref{EulerRDE}}^{2}=C_{\ref{EulerRDE}}^{2}\left(
p,N\right) $. As in part (a) we can find Lipschitz paths $\bar{x}^{n}$ such
that%
\begin{equation*}
S_{N}\left( \bar{x}^{n}\right) \rightarrow \mathbf{\bar{x}}
\end{equation*}%
uniformly on $\left[ 0,1\right] $, such that%
\begin{equation*}
\left\Vert S_{N}\left( \bar{x}^{n}\right) \right\Vert _{1/p\text{-H\"{o}lder}%
}\leq 3C_{\ref{EulerRDE}}^{2}M\equiv M_{1}
\end{equation*}%
Note that, by projection, $S_{\left[ p\right] }\left( \bar{x}^{n}\right)
\rightarrow \mathbf{x}$ uniformly on $\left[ 0,1\right] $ with uniform
homogenous $1/p$-H\"{o}lder bounds. As before, the Universal Limit Theorem
implies that%
\begin{equation*}
\pi \left( 0,y_{0};\bar{x}^{n}\right) \rightarrow \pi \left( 0,y_{0};\mathbf{%
x}\right) \text{ \ uniformly on }\left[ 0,1\right]
\end{equation*}%
while Davie's lemma implies the existence of $C_{\ref{EulerRDE}}^{3}$, de
facto $C_{\ref{tfa}}=C_{\ref{tfa}}\left( M_{1}\right) $ from lemma \ref{tfa}%
, such that, uniformly over $n$, and for all $0\leq s<t\leq 1$,%
\begin{equation*}
\left\vert \pi \left( 0,y_{0};\bar{x}^{n}\right) _{s,t}-I\left[ \pi \left(
0,y_{0};\bar{x}^{n}\right) _{s},N,S_{N}\left( \bar{x}^{n}\right) _{s,t}%
\right] \right\vert \leq C_{\ref{EulerRDE}}^{3}\omega \left( s,t\right)
^{\theta }
\end{equation*}%
with $\theta =\left( N+1\right) /p$. By continuity of the map 
\begin{equation*}
\left( z,\mathbf{g}\right) \in \mathbb{R}^{e}\times G^{N}\left( \mathbb{R}%
^{d}\right) \mapsto I\left[ z,N,\mathbf{g}\right] \in \mathbb{R}^{e},
\end{equation*}%
we can send $n\rightarrow \infty $ to obtain%
\begin{equation*}
\left\vert \pi \left( 0,y_{0};\mathbf{x}\right) _{s,t}-I\left[ \pi \left(
0,y_{0};\mathbf{x}\right) _{s},N,\mathbf{x}_{s,t}\right] \right\vert \leq C_{%
\ref{EulerRDE}}^{3}\omega \left( s,t\right) ^{\theta }\text{.}
\end{equation*}%
Finally, as above,%
\begin{equation*}
\ln C_{\ref{EulerRDE}}^{3}=O\left( \left( \ln M_{1}\right) ^{2}\right)
=O\left( \left( \ln M\right) ^{2}\right) \text{.}
\end{equation*}
\end{proof}

\section{Asymptotic Expansions for RDE Flows}

We now consider RDEs driven by a random geometric $p$-rough path $\mathbf{x=x%
}\left( \omega \right) $ defined on some complete probability space $\left(
\Omega ,\mathcal{F},\mathbb{P}\right) .$ We shall assume that the r.v. $%
\left\Vert \mathbf{x}\right\Vert _{1/p\text{-H\"{o}lder;}\left[ 0,1\right] }$
has Gauss tails since this is the case for all examples we have in mind:
Enhanced Brownian motion $\mathbf{B}$, see \cite{FV03, FLS}, Enhanced
Fractional Brownian Motion $\mathbf{B}^{H}$ and other Enhanced Gaussian
processes \cite{CQ, FV04} and Enhanced Markov processes with uniformly
elliptic generated in divergence form \cite{Le}.\ However, the proof of the
following theorem will make clear that the method works whenever the
real-valued r.v. $\left\Vert \mathbf{x}\right\Vert _{1/p\text{-H\"{o}lder;}%
\left[ 0,1\right] }$ has some exponential tail decay.

\begin{theorem}
\label{azencott}Let $p\in \left( 2,3\right) $, $N\geq \left[ p\right] +1$
and consider the random RDE solution $\pi \left( 0,y_{0},\mathbf{x}\right) $
driven by the random geometric $p$-rough path $\mathbf{x}$ (along fixed $%
\mathrm{Lip}^{N}$ vector fields $V_{1},...,V_{d}$), assuming that $%
\left\Vert \mathbf{x}\right\Vert _{1/p\text{-H\"{o}lder;}\left[ 0,1\right] }$
has Gauss tails,%
\begin{equation*}
\exists \alpha >0:\mathbb{E}\left[ \exp \left( \alpha \,\left\Vert \mathbf{x}%
\right\Vert _{1/p\text{-H\"{o}lder;}\left[ 0,1\right] }^{2}\right) \right]
<\infty .
\end{equation*}%
Then 
\begin{equation*}
\left\vert \pi \left( 0,y_{0},\mathbf{x}\right) \right\vert _{1/p\text{-H%
\"{o}lder;}\left[ 0,1\right] }\in L^{q}\left( \Omega \right) \text{ for all }%
q\in \lbrack 1,\infty ).
\end{equation*}%
Moreover, the remainder of the step-$N$ Euler approximation is bounded in
probability. More precisely, there is a constant $C_{\ref{azencott}}$
dependent on \thinspace \thinspace $\alpha ,p,N,y_{0}$ and $V_{1},...,V_{d}$
such that for $R\geq 1$ and all $t\in (0,1],$ 
\begin{equation*}
\mathbb{P}\left( \sup_{0\leq s\leq t}\left\vert \pi \left( 0,y_{0},\mathbf{x}%
\right) _{0,s}-I\left[ y_{0},N,S_{N}\left( \mathbf{x}\right) _{0,s}\right]
\right\vert >Rt^{\frac{N+1}{p}}\right) \leq C_{\ref{azencott}}\exp \left(
-e^{\left( \ln R\right) ^{1/2}/C_{\ref{azencott}}}\right) .
\end{equation*}%
In particular, the l.h.s. tends to zero uniformly over $t\in (0,1]$ as $%
R\rightarrow \infty $ and the convergence is faster than any power of $1/R$.
\end{theorem}

\begin{proof}
By assumption, $M_{1}=\max \left\{ \left\Vert \mathbf{x}\right\Vert _{1/p%
\text{-H\"{o}lder;}\left[ 0,1\right] },1\right\} $ has Gauss tails and $%
\mathbb{E}\left[ \exp \left( \alpha \,M_{1}^{2}\right) \right] <\infty $.
From Theorem \ref{EulerRDE},%
\begin{equation*}
\sup_{0\leq s\leq t}\left\vert \pi \left( 0,y_{0},\mathbf{x}\right)
_{0,s}-Iy_{0},N,S_{N}\left( \mathbf{x}\right) _{0,s}\right\vert \leq C_{%
\text{\ref{EulerRDE}}}^{\prime }\left( M_{1}\right) t^{\left( N+1\right) /p}.
\end{equation*}%
where 
\begin{equation*}
C_{\text{\ref{EulerRDE}}}^{\prime }\left( M_{1}\right) \leq C_{\text{\ref%
{azencott}}}^{1}e^{\left( C_{\text{\ref{azencott}}}^{1}\ln M_{1}\right)
^{2}}.
\end{equation*}%
Therefore, 
\begin{eqnarray*}
&&\mathbb{P}\left[ \sup_{0\leq s\leq t}\left\vert \pi \left( 0,y_{0},\mathbf{%
x}\right) _{0,s}-I^{y_{0},N,S_{N}\left( \mathbf{x}\right) _{0,s}}\right\vert
>Rt^{\frac{N+1}{p}}\right] \\
&\leq &\mathbb{P}\left[ C_{\text{\ref{azencott}}}^{1}e^{\left( C_{\text{\ref%
{azencott}}}^{1}\ln M_{1}\right) ^{2}}>R\right] \\
&=&\mathbb{P}\left[ M>\exp \left( \frac{\sqrt{\ln \left( R/C_{\text{\ref%
{azencott}}}^{1}\right) }}{C_{\text{\ref{azencott}}}^{1}}\right) \right] \\
&\leq &\mathbb{E}\left[ \exp \left( \alpha \,M^{2}\right) \right] \exp \left[
-\alpha \exp \left( \frac{2\sqrt{\ln \left( R/C_{\text{\ref{azencott}}%
}^{1}\right) }}{C_{\text{\ref{azencott}}}^{1}}\right) \right] \\
&\leq &C_{\text{\ref{azencott}}}^{2}\exp \left( -e^{\left( \ln \left(
R\right) \right) ^{1/2}/C_{\text{\ref{azencott}}}^{2}}\right) ,
\end{eqnarray*}%
where the last estimate is valid for every $R\geq 1$ by choosing $C_{\text{%
\ref{azencott}}}^{2}$ sufficiently large.
\end{proof}

\begin{remark}
The same result holds if we replace $I\left[ y_{0},N,S_{N}\left( \mathbf{x}%
\right) _{0,s}\right] $ by $\pi \left( 0,y_{0},x^{0,s}\right) _{0,s},$ where 
$x^{0,s}$ is a geodesic associated to the element $S_{N}\left( \mathbf{x}%
\right) _{0,s}$ of $G^{N}\left( \mathbb{R}^{d}\right) $.
\end{remark}

\begin{remark}
Even in the case of Enhanced Brownian motion, $\mathbf{x=B}$, probability
estimates of the unrestricted event%
\begin{equation*}
\left\{ \sup_{0\leq s\leq t}\left\vert \pi \left( 0,y_{0},\mathbf{x}\right)
_{0,s}-I\left[ y_{0},N,S_{N}\left( \mathbf{x}\right) _{0,s}\right]
\right\vert >Rt^{\frac{N+1}{p}}\right\}
\end{equation*}
valid for all $\left( t,R\right) \in (0,1]\times $ $[1,\infty )$ appear
novel compared to the results given in \cite{Az, BA, Ca}.
\end{remark}

The perhaps strongest estimate that has been extracted from Azencott's work
in this context (see \cite[p 235]{Ca}) is the following: in our notation
(recall that $\pi \left( 0,y_{0},\mathbf{B}\right) $ solves a Stratonovich
stochastic differential equation): $\exists a,c>0:\forall R\geq 0:$%
\begin{equation}
\underset{t\rightarrow 0}{\overline{\lim }}\mathbb{P}\left( \sup_{0\leq
s\leq t}\left\vert \pi \left( 0,y_{0},\mathbf{B}\right)
_{0,s}-I^{y_{0},N,S_{N}\left( \mathbf{B}\right) _{0,s}}\right\vert >Rt^{%
\frac{N+1}{2}}\right) \leq ce^{-\frac{R^{a}}{c}}.
\label{AzencottCastellEstimate}
\end{equation}%
(Note that the exponent of $t$ is$\left( N+1\right) /2$ in contrast to $%
\left( N+1\right) /p$ in Theorem \ref{azencott}). We now show how (\ref%
{AzencottCastellEstimate}) can be deduced from our general results.

\begin{proposition}
\label{CastellCor}We keep all assumptions of the preceding theorem but drive
the RDE with Enhanced Brownian motion $\mathbf{x}=\mathbf{B}$. Then (\ref%
{AzencottCastellEstimate}) holds with $a=2/\left( N+1\right) $ and $c=C_{%
\text{\ref{CastellCor}}}$ depending on \thinspace $N,y_{0}$ and $%
V_{1},...,V_{d}$.
\end{proposition}

\begin{proof}
Choose $p=p\left( N\right) $ s.t.%
\begin{equation*}
2<p<2\frac{N+2}{N+1}\leq 3
\end{equation*}%
Then there exists $\alpha =\alpha \left( p\right) >0$ s.t. $\mathbb{E}\left[
\exp \left( \alpha \,M^{2}\right) \right] <\infty $ where $M=\max \left\{
\left\Vert \mathbf{B}\right\Vert _{1/p\text{-H\"{o}lder;}\left[ 0,1\right]
},1\right\} $. We set $\mathbf{\bar{B}\equiv }$ $S_{N+1}\left( \mathbf{B}%
\right) $, well-defined by Proposition \ref{th1}. Then

\begin{eqnarray*}
\left\vert \pi \left( 0,y_{0},\mathbf{B}\right)
_{0,s}-I^{y_{0},N,S_{N}\left( \mathbf{B}\right) _{0,s}}\right\vert &\leq
&\left\vert \pi \left( 0,y_{0},\mathbf{B}\right) _{0,s}-I^{y_{0},N+1,\mathbf{%
\bar{B}}_{0,s}}\right\vert \\
&&+\left\vert I^{y_{0},N+1,\mathbf{\bar{B}}_{0,s}}-I^{y_{0},N,S_{N}\left( 
\mathbf{x}\right) _{0,s}}\right\vert .
\end{eqnarray*}%
From Theorem \ref{EulerRDE} and Proposition \ref{th1},%
\begin{equation*}
\left\vert \pi \left( 0,y_{0},\mathbf{x}\right) _{0,s}-I^{y_{0},N+1,\mathbf{%
\bar{B}}_{0,s}}\right\vert \leq C_{\text{\ref{CastellCor}}}^{1}s^{\left(
N+2\right) /p}e^{\left( C_{\text{\ref{CastellCor}}}^{1}\ln M\right) ^{2}}.
\end{equation*}%
On the other hand%
\begin{eqnarray*}
\left\vert I^{y_{0},N+1,S_{N+1}\left( \mathbf{B}\right)
_{0,s}}-I^{y_{0},N,S_{N}\left( \mathbf{x}\right) _{0,s}}\right\vert &=&|\sum 
_{\substack{ i_{1},...,i_{N+1}  \\ \in \left\{ 1,...,d\right\} }}%
V_{i_{1}}\cdots V_{i_{N+1}}H\left( y_{0}\right) \mathbf{\bar{B}}%
_{0,s}^{N+1,i_{1},\cdots ,i_{N+1}}|\text{ } \\
&\leq &C_{\text{\ref{CastellCor}}}^{2}\sum_{i_{1},...,i_{N+1}\in \left\{
1,...,d\right\} }\left\vert \mathbf{\bar{B}}_{0,s}^{N+1,i_{1},\cdots
,i_{N+1}}\right\vert \\
&\leq &C_{\text{\ref{CastellCor}}}^{3}\left\Vert \mathbf{\bar{B}}%
_{0,s}\right\Vert ^{N+1}.
\end{eqnarray*}%
Trivially, $\sup_{0\leq s\leq t}s^{\left( N+2\right) /p}=t^{\left(
N+2\right) /p}$ and we are led to%
\begin{eqnarray*}
&&\mathbb{P}\left( \sup_{s\in \left[ 0,t\right] }\left\vert \pi \left(
0,y_{0},\mathbf{x}\right) _{0,s}-I^{y_{0},N,S_{N}\left( \mathbf{x}\right)
_{0,s}}\right\vert \geq Rt^{\frac{N+1}{2}}\right) \\
&\leq &\mathbb{P}\left( t^{\frac{N+2}{p}}C_{\text{\ref{CastellCor}}%
}^{1}e^{\left( C_{\text{\ref{CastellCor}}}^{1}\ln M\right) ^{2}}+C_{\text{%
\ref{CastellCor}}}^{3}\sup_{s\in \left[ 0,t\right] }\left\Vert \mathbf{\bar{B%
}}_{0,s}\right\Vert ^{N+1}\geq Rt^{\frac{N+1}{2}}\right) \\
&\leq &\left( 1\right) +(2)
\end{eqnarray*}%
where with $\nu =\left( N+2\right) /p-\left( N+1\right) /2>0,$%
\begin{equation*}
\left( 1\right) =\mathbb{P}\left( t^{\nu }C_{\text{\ref{CastellCor}}%
}^{1}e^{\left( C_{\text{\ref{CastellCor}}}^{1}\ln M\right) ^{2}}\geq \frac{R%
}{2}\right)
\end{equation*}%
and%
\begin{equation*}
\left( 2\right) =\mathbb{P}\left( C_{\text{\ref{CastellCor}}}^{3}\left(
t^{-1/2}\sup_{s\in \left[ 0,t\right] }\left\Vert \mathbf{\bar{B}}%
_{0,s}\right\Vert \right) ^{N+1}\geq \frac{R}{2}\right) .
\end{equation*}%
Gauss tails of $M$ are more than enough to asset that $(1)$ tends to zero as 
$t\rightarrow 0$. As for (2), Brownian scaling shows that $\left( 2\right) $
is in fact independent of $t$ and hence equal to%
\begin{equation*}
\mathbb{P}\left( C_{\text{\ref{CastellCor}}}^{3}\left( \sup_{s\in \left[ 0,1%
\right] }\left\Vert \mathbf{\bar{B}}_{0,s}\right\Vert \right) ^{N+1}\geq 
\frac{R}{2}\right) \leq \mathbb{P}\left( C_{\text{\ref{CastellCor}}%
}^{4}M^{N+1}\geq \frac{R}{2}\right)
\end{equation*}%
since $\sup_{s\in \left[ 0,1\right] }\left\Vert \mathbf{\bar{B}}%
_{0,s}\right\Vert \leq \left\Vert \mathbf{\bar{B}}\right\Vert _{1/p\text{-H%
\"{o}lder;}\left[ 0,1\right] }\leq $ $C_{\text{\ref{CastellCor}}%
}^{5}\left\Vert \mathbf{B}\right\Vert _{1/p\text{-H\"{o}lder;}\left[ 0,1%
\right] }\leq C_{\text{\ref{CastellCor}}}^{5}M$. Gauss tails of $M$ now
easly the claimed tail decay.
\end{proof}

\begin{remark}
Estimate (\ref{AzencottCastellEstimate}) remains valid when $I\left[
y_{0},N,S_{N}\left( \mathbf{B}\right) _{0,s}\right] $ is replaced by $\pi
\left( 0,y_{0};x^{0,s}\right) $ where $x^{0,s}$ is a geodesic associated to $%
S_{N}\left( \mathbf{x}\right) _{0,s}.$
\end{remark}

\begin{remark}
Note that $a=2/\left( N+1\right) $ is what one expects from integrability of
the $\left( N+1\right) ^{th}$ multiple Wiener-It\^{o} integral which
dominates the remainder as $t\rightarrow 0$. Theorem \ref{azencott} gives an
estimate valid uniformly in $t\in (0,1]$ and we would not expect an
exponential tail. In fact, the given estimate implies a tail decay which is
about as good as one can hope in absence of exponential decay and there
seems little room for improvement.
\end{remark}

\begin{remark}
Only in absence of a drift vector field $V_{0}$ does our step-$N$ Euler
approximation $I\left[ y_{0},N,S_{N}\left( \mathbf{x}\right) _{0,s}\right] $
coincide \textit{precisely} with the approximation of \cite[prop 4.3]{Az}
and \cite[p 235]{Ca}. When $V_{0}\neq 0$, the methodology and results of $%
\left( p,q\right) $-rough paths \cite{LV} could be used to produce the same
approximations as those in the above cited references. If one accepts a few
additional terms in the step-$N$ approximation, it may be simplest to deal
with $V_{0}\neq 0$ via an RDE\ driven by the canonically defined time-space
rough path.
\end{remark}

\begin{remark}
Proposition \ref{CastellCor} is readily adapted to EFBM $\mathbf{B}^{H}$
with $H>1/4$ and gives%
\begin{equation*}
\underset{t\rightarrow 0}{\overline{\lim }}\mathbb{P}\left( \sup_{0\leq
s\leq t}\left\vert \pi \left( 0,y_{0},\mathbf{B}^{H}\right)
_{0,s}-I^{y_{0},N,S_{N}\left( \mathbf{B}^{H}\right) _{0,s}}\right\vert
>Rt^{H\left( N+1\right) }\right) \leq ce^{-\frac{R^{a}}{c}}.
\end{equation*}
\end{remark}

\begin{remark}
In comparison to the Taylor expansion of Azencott, the approximations by Ben
Arous and Castell respect the geometry of the problem. This is also the case
for our geodesic approximations although the efficient approximations of
sub-Riemannian geodesics remains a numerical challenge.
\end{remark}

\section{$L^{q}$ Convergence in the Universal Limit Theorem}

We now give a criterion for $L^{q}$-convergence in the Universal Limit
Theorem.

\begin{proposition}
\bigskip Assume that for a random sequence of rough path%
\begin{equation*}
\sup_{n}\left\Vert \mathbf{x}_{n}\right\Vert _{1/p\text{-H\"{o}lder;}\left[
0,1\right] }\equiv M
\end{equation*}%
has a tail Gaus tail. Assume that $d_{1/p\text{-H\"{o}lder}}\left( \mathbf{x}%
_{n},\mathbf{x}\right) \rightarrow 0$ in probabililty. Then%
\begin{equation*}
\left\vert \pi \left( 0,y_{0},\mathbf{x}_{n}\right) -\pi \left( 0,y_{0},%
\mathbf{x}\right) \right\vert _{1/p\text{-H\"{o}lder;}\left[ 0,1\right]
}\rightarrow 0\text{ in }L^{q}\text{ }\forall q\in \lbrack 1,\infty )\text{.}
\end{equation*}
\end{proposition}

\begin{proof}
By the universal limit theorem,%
\begin{equation*}
Z_{n}\equiv \left\vert \pi \left( 0,y_{0},\mathbf{x}_{n}\right) -\pi \left(
0,y_{0},\mathbf{x}\right) \right\vert _{1/p\text{-H\"{o}lder;}\left[ 0,1%
\right] }\rightarrow 0\text{ in probability.}
\end{equation*}%
At least along a subsequence, $Z_{n_{k}}\rightarrow 0$ a.s. as $k\rightarrow
\infty $ and%
\begin{equation*}
\left\vert \pi \left( 0,y_{0},\mathbf{x}\right) \right\vert _{1/p\text{-H%
\"{o}lder;}\left[ 0,1\right] }=\lim_{k\rightarrow \infty }\left\vert \pi
\left( 0,y_{0},\mathbf{x}_{n_{k}}\right) \right\vert _{1/p\text{-H\"{o}lder;}%
\left[ 0,1\right] }
\end{equation*}%
and therefore 
\begin{equation*}
Z_{n}\leq 2\sup_{n}\left\vert \pi \left( 0,y_{0},\mathbf{x}_{n}\right)
\right\vert _{1/p\text{-H\"{o}lder;}\left[ 0,1\right] }\leq 2C
\end{equation*}%
where $C=C\left( M\right) \in $ $L^{q}$ $\forall q\in \lbrack 1,\infty )$
using Theorem \ref{azencott}. To obtain $L^{q}$ convergence we have to check
that $\left\{ Z_{n}^{q}:n\geq 1\right\} $ is uniformly integrable. But this
follows immediately from $L^{r}$-boundedness, $r>1$, indeed%
\begin{equation*}
\mathbb{E}\left( \left\vert Z_{n}^{q}\right\vert ^{r}\right) \leq 2^{qr}%
\mathbb{E}\left( C\left( M\right) ^{qr}\right) <\infty .
\end{equation*}
\end{proof}

\begin{example}
Let $\mathbf{x=B}$ be Enhanced Brownian motion over $\mathbb{R}^{d}$ with%
\begin{equation*}
\left\Vert \mathbf{x}\right\Vert _{1/p\text{-H\"{o}lder;}\left[ 0,1\right]
}<\infty \text{ where }p\in \left( 2,3\right) .
\end{equation*}%
Let $\left( D_{n}\right) $ be a nested family of dissections of $\left[ 0,1%
\right] $ with mesh $\left\vert D_{n}\right\vert \rightarrow 0$ and define
the lifted piecewise linear approximation%
\begin{equation*}
\mathbf{x}_{n}=S_{2}\left( \mathbb{E}\left[ B\mathbf{|}\sigma \left(
B_{t}:t\in D_{n}\right) \right] \right) .
\end{equation*}%
Then observe that the area process of $\mathbf{x}_{n}$ is obtained by the L%
\'{e}vy area process condtioned on $\sigma \left( B_{t}:t\in D_{n}\right) $.
Then 
\begin{equation*}
\sup_{n}\left\Vert \mathbf{x}_{n}\right\Vert _{1/p\text{-H\"{o}lder;}\left[
0,1\right] }\leq \sup_{n}\mathbb{E}\left[ \left\Vert \mathbf{x}\right\Vert
_{1/p\text{-H\"{o}lder;}\left[ 0,1\right] }|\sigma \left( B_{t}:t\in
D_{n}\right) \right] \equiv \sup_{n}M\left( n\right) \equiv M<\infty .
\end{equation*}%
Note that 
\begin{eqnarray*}
M\left( n\right) &\equiv &\mathbb{E}\left[ \left\Vert \mathbf{x}\right\Vert
_{1/p\text{-H\"{o}lder;}\left[ 0,1\right] }|\sigma \left( B_{t}:t\in
D_{n}\right) \right] \\
&\leq &\mathbb{E}\left[ \left\Vert \mathbf{x}\right\Vert _{1/p\text{-H\"{o}%
lder;}\left[ 0,1\right] }^{q}|\sigma \left( B_{t}:t\in D_{n}\right) \right]
^{1/q}
\end{eqnarray*}%
and from Doob's maximal inequality,%
\begin{eqnarray*}
\mathbb{E}\left\vert \sup_{n}\left\Vert \mathbf{x}_{n}\right\Vert _{1/p\text{%
-H\"{o}lder;}\left[ 0,1\right] }^{q}\right\vert &\leq &\left( \frac{q}{q-1}%
\right) ^{q}\sup_{n}\mathbb{E}\left( M\left( n\right) ^{q}\right) \\
&\leq &\left( \frac{q}{q-1}\right) ^{q}\mathbb{E}\left[ \left\Vert \mathbf{x}%
\right\Vert _{1/p\text{-H\"{o}lder;}\left[ 0,1\right] }^{q}\right] .
\end{eqnarray*}%
Noting that $\left( q/\left( q-1\right) \right) ^{q}$ stays bounded as $%
q\rightarrow \infty $ (in fact, it converges to $e$), we conclude that Gauss
tails of $M_{1}$ imply Gauss tails of $M$. Finally, fix $p\in \left(
2,3\right) $ and repeat the above argument for $\tilde{p}\in \left(
2,p\right) $ replacing $p$. By interpolation%
\begin{equation*}
d_{1/p\text{-H\"{o}lder}}\left( \mathbf{x}_{n},\mathbf{x}\right) \rightarrow
0\text{ a.s.}
\end{equation*}%
and uniform Gauss tails follow from the trivial estimate $\left\Vert \mathbf{%
x}_{n}\right\Vert _{1/p\text{-H\"{o}lder;}\left[ 0,1\right] }\leq \left\Vert 
\mathbf{x}_{n}\right\Vert _{1/\tilde{p}\text{-H\"{o}lder;}\left[ 0,1\right]
} $.
\end{example}

\begin{remark}
The same argument works for other Enhanced Gaussian processes with
martingale approximations such as fractional Brownian motion, see \cite{FV04}
for the case $H>1/3$.
\end{remark}

\end{document}